%%%%%%%              THIS IS AN AMSTEX-FILE        %%%%%%%
%%%%%%%                                            %%%%%%%
\documentstyle{amsppt}
\magnification=1200
\baselineskip=18pt
\nologo
\TagsOnRight
%\hsize=15 true cm
\document
\define \ts{\thinspace}
\define \di{\partial}

\define \ot{\otimes}

\define \gl{{\frak {gl}}}
\define \U{\operatorname {U}}

\define \Mat{{\operatorname {Mat}}}

\define \C{\Bbb C}
\define \PD{{\Cal {PD}}}
\define \Proof{\noindent {\bf Proof. }}

\heading{\bf A REMARK ON THE HIGHER CAPELLI IDENTITIES}
\endheading
\bigskip
\bigskip
\heading{Alexander Molev}\endheading
\bigskip
\bigskip
\bigskip
\heading{\sl Dedicated to Jimmy}\endheading
\bigskip
\bigskip
\bigskip
\bigskip
\noindent
Centre for Mathematics and its Applications,\newline
Australian National University,\newline
Canberra, ACT 0200, Australia\newline
(e-mail: molev\@pell.anu.edu.au)
\bigskip
\bigskip
\noindent
{\bf February 1996}
\bigskip
\bigskip
\noindent
{\bf Abstract}\newline
A simple proof of the higher Capelli identities is given.

\newpage

The aim of this note is to give a simple proof of the remarkable
generalizations of the classical Capelli identity discovered by
A. Okounkov [O1] and also proved in different ways
by him [O2] and M. Nazarov [N]. 
Following [O2], we use some properties of the Jucys--Murphy
elements in the group algebra for the symmetric group,
together with the branching properties of the Young basis.
The difference between our proof
and that of [O2] is that in our approach
we do not need to use the Wick formula and the Olshanski\u\i\
special symmetrization map.
\medskip

Denote by $\PD$ the algebra of polynomial coefficient differential
operators in $mn$ variables $x_{ai}$, where $a=1,\dots,m$ and
$i=1,\dots,n$. Consider the representation of the Lie
algebra $\gl(m)$ by differential operators defined on
the standard generators $E_{ab}$ as follows:
$$
E_{ab}=\sum_{i=1}^n x_{ai}\di_{bi},\tag 1
$$
where $\di_{ai}:=\di/\di x_{ai}$. Let $E$, $X$, $D$ denote the formal
matrices with the entries $E_{ab}$, $x_{ai}$, $\di_{ai}$,
respectively. Then (1) can be
written in a matrix form as follows:
$$
E=XD', \tag 2
$$
where $D'$ is the matrix transposed to $D$.

As in [O2], for a partition
$\lambda$ and a standard $\lambda$-tableau $T$ we denote by $v_T$
the corresponding vector of the Young orthonormal basis in the
irreducible representation $V^{\lambda}$ of the symmetric group $S_k$,
$k=|\lambda|$. We let $c_T(r)=j-i$ if the cell $(i,j)\in\lambda$
is occupied by the entry $r$ of the tableau $T$. Introduce the
matrix element
$$
\Psi_{TT'}=\sum_{s\in S_k}(s\cdot v_T,v_{T'})\cdot s^{-1}\in\C[S_k].\tag 3
$$
The symmetric group $S_k$ acts in a natural way in the tensor space
$(\C^m)^{\ot k}$,
so that we can identify permutations from $S_k$ with
elements of the algebra
$$
\underbrace{\Mat_{mm}\ot\cdots\ot\Mat_{mm}}_k,\tag 4
$$
where by $\Mat_{pq}$ we denote the space of $p\times q$-matrices;
for $p=q$ we also regard it as an algebra. By $P_{ij}$ we denote
the element of (4) corresponding to the transposition $(ij)\in S_k$.

We regard the tensor product $A\ot B\ot\cdots\ot C$
of $k$ matrices $A$, $B$, $\dots$, $C$ 
of size $p\times q$ 
with entries from an algebra $\Cal A$ as an element
$$
\sum A_{a_1i_1}B_{a_2i_2}\cdots C_{a_ki_k}
\ot e_{a_1i_1}\ot e_{a_2i_2}\ot\cdots\ot e_{a_ki_k}\in\Cal A\ot
(\Mat_{pq})^{\ot k},
$$
where the $e_{ai}$ are the standard matrix units.

\bigskip
\proclaim{{\bf Theorem} {\rm [O2, N]}} Let $T$ and $T'$ be
two standard tableaux of the same shape. Then
$$
(E-c_T(1))\ot\cdots\ot(E-c_T(k))\cdot\Psi_{TT'}=
X^{\ot k}\cdot(D')^{\ot k} \cdot\Psi_{TT'}.\tag 5
$$
\endproclaim

\Proof We use induction on $k$. Denote by $U$ the tableau
obtained from $T$ by removing the cell with the entry $k$.
Using the branching property of the Young basis $\{v_T\}$ one
can easily check that
$$
\Psi_{TT'}=\text{\rm const}\cdot\Psi_{UU}\Psi_{TT'},
$$
where `const' is a nonzero constant (more precisely, $\text{const}=
\text{dim}\ts \mu/(k-1)!$ where $\mu$ is the shape of $U$
and $\text{dim}\ts \mu=\text{dim}\ts V^{\mu}$).

So, we can rewrite the
left hand side of (5) as follows:
$$
\text{const}\cdot(E-c_T(1))\ot\cdots\ot(E-c_T(k-1))\cdot\Psi_{UU}
\ot(E-c_T(k))\cdot\Psi_{TT'}.
$$
By the induction hypothesis, this equals
$$
\text{\rm const}\cdot X^{\ot k-1}\cdot(D')^{\ot k-1} \cdot\Psi_{UU}
\ot(XD'-c_T(k))\cdot\Psi_{TT'}
$$
$$
=X^{\ot k-1}\cdot(D')^{\ot k-1}
\ot(XD'-c_T(k))\cdot\Psi_{TT'}
$$
$$
=\bigl(\sum x_{a_1i_1}\cdots x_{a_{k-1}i_{k-1}}
\di_{b_1i_1}\cdots \di_{b_{k-1}i_{k-1}}
(x_{a_ki_k}\di_{b_ki_k}-\frac{\delta_{a_kb_k}}{n} c_T(k))
$$
$$
{}\ot e_{a_1b_1}\ot\cdots\ot e_{a_kb_k}\bigr)\cdot\Psi_{TT'}.
$$
Now we transform this expression using the relations 
$\di_{bj}x_{ai}=x_{ai}\di_{bj}+\delta_{ab}\delta_{ij}$
to obtain
$$
\gather
\left(\sum x_{a_1i_1}\cdots x_{a_{k}i_{k}}
\di_{b_1i_1}\cdots \di_{b_{k}i_{k}}
\ot e_{a_1b_1}\ot\cdots\ot e_{a_kb_k}\right)\cdot\Psi_{TT'}\\
+\left(\sum x_{a_1i_1}\cdots x_{a_{k-1}i_{k-1}}
\di_{b_1i_1}\cdots \di_{b_{k-1}i_{k-1}}
\ot e_{a_1b_1}\ot\cdots\ot e_{a_{k-1}b_{k-1}}\ot 1\right)\\
\times(P_{1k}+\cdots+P_{k-1,k}-c_T(k))\cdot\Psi_{TT'}.\tag 6
\endgather
$$
Note that $P_{1k}+\cdots+P_{k-1,k}$ is the image of the Jucys--Murphy
element $(1k)+\cdots+(k-1,k)$. It has the property
$$
((1k)+\cdots+(k-1,k))\cdot\Psi_{TT'}=c_T(k)\cdot\Psi_{TT'},
$$
which was also used in [O2] and can be easily derived from the following
formula due to Jucys and Murphy:
$$
((1k)+\cdots+(k-1,k))\cdot v_T=c_T(k)\cdot v_T.
$$
This proves that the second summand in (6)
is zero, which completes the proof of
Theorem. 
\medskip

It was shown in [O2] that 
taking trace in both sides of (5) over all tensor factors $\Mat_{mm}$, 
one obtains the 
following `higher Capelli identities' proved in [O1] and [N].
\bigskip
\proclaim
{\bf Corollary} For any standard tableau $T$ of shape $\lambda$
one has
$$
\text{\rm tr}\ts (E-c_T(1))\ot\cdots\ot(E-c_T(k))\cdot\Psi_{TT}=
\frac1{\text{\rm dim}\ts \lambda}\ts
\text{\rm tr}\ts X^{\ot k}\cdot(D')^{\ot k}\cdot\chi^{\lambda},
\tag 7
$$
where
$$
\chi^{\lambda}=\sum_{s\in S_k}\chi^{\lambda}(s)\cdot s\in\C[S_k]
$$
is the character of $V^{\lambda}$.
\endproclaim

In particular, in the case $\lambda=(1^k)$ one obtains
the classical Capelli identity.

The left hand side of (7) is an element of the center of the universal
enveloping algebra $\U(\gl(m))$. Due to (7), it depends only on the
partition $\lambda$ and does not depend on the tableau $T$. Moreover,
the set of these elements, where $\lambda$ runs over the partitions
with length $\leq m$ forms a basis in the center. Their
eigenvalues in highest weight representations of $\gl(m)$ are
the `shifted Schur polynomials'; see [O1, O2, OO, N] for details.

\bigskip
\noindent
{\bf References}
\bigskip

\itemitem{[N]}
{M.~Nazarov}, {Yangians and Capelli identities}, preprint, 1995;
q-alg/9601027.

\itemitem{[O1]}
A.~Okounkov,
{Quantum immanants and higher Capelli identities},
preprint, 1995;\newline
q-alg/9602028.

\itemitem{[O2]}
A.~Okounkov,
{Young basis, Wick formula and higher Capelli identities},
preprint, 1995; \newline
q-alg/9602027.

\itemitem{[OO]}
A.~Okounkov and G.~Olshanski\u\i,
{Shifted Schur functions},
preprint, 1995.

\enddocument